\documentclass[11pt]{article}
\usepackage{lscape}
\usepackage{afterpage}
\usepackage{pstricks}
\usepackage{pst-plot}
\usepackage{longtable}
\usepackage{dcolumn}
\usepackage{pst-node}
\usepackage[dvips]{graphicx}
\usepackage{amsmath}
\usepackage{amssymb}
\usepackage{amsthm}
\usepackage{multirow}
\usepackage{caption}
\usepackage{colortab}
\usepackage{color}
\usepackage{array}
\usepackage{slashbox}
\usepackage{colortbl}
\usepackage{shadbox}
\usepackage{subfigure}
\usepackage{textcomp}
\usepackage{setspace}
\usepackage{pstricks, pst-node}
\psset{arrows=->, labelsep=3pt, mnode=circle}
\usepackage{eurosym}

\newfont{\rams}{msbm10 scaled\magstep1}
\newcommand{\rea}{\mbox{\rams \symbol{'122}}}

\newcommand{\rio}{\mathbb{R}}

\newenvironment{resumeT}{\begin{list}{}{\setlength{\rightmargin}{\leftmargin}}\item[]
{\centering {\bf \it~~~}
\par}\item[]\ignorespaces}{\unskip\end{list}}

\newtheorem{theorem}{Theorem}[section]
\newtheorem{note}{Remark}[section]

\begin{document}

\title{\Large
    \textbf{The SMAA-PROMETHEE methods}
  }  
  
\author{ \hspace{0,1cm} Salvatore Corrente \thanks{Department of Economics and Business, University of Catania, Corso Italia 55, 95129  Catania, Italy, e-mails: \texttt{salvatore.corrente\string@unict.it}, \texttt{salgreco\string@unict.it}},
Jos\'e Rui Figueira \footnote{CEG-IST, Instituto Superior T\'{e}cnico, Technical University of Lisbon, Av. Rovisco Pais, 1049-001 Lisboa, Portugal, e-mail: \texttt{figueira@ist.utl.pt}}, \hspace{0,1cm} Salvatore Greco $^*$}


\date{}
\maketitle 


\addcontentsline{toc}{section}{Abstract}


\vspace{-1,5cm}

\begin{resumeT}

\textbf{Abstract:} \noindent PROMETHEE methods are widely used in Multiple Criteria Decision Aiding (MCDA) to deal with real decision making problems. A crucial aspect of the classical PROMETHEE methods is the choice of criteria weights. In this paper, we propose to apply the Stochastic Multiobjective Acceptability Analysis (SMAA) to the classical PROMETHEE methods and to the bipolar PROMETHEE methods in order to explore the whole set of weights compatible with some preference information provided by the Decision Maker (DM). A didactic example describes the application of the presented methodology to a student evaluation problem.

\end{resumeT}

 \pagenumbering{arabic}

\section{Introduction}
PROMETHEE is a well-known family of Multiple Criteria Decision Aiding (MCDA) outranking methods (see \cite{FigGreEhr} for a survey on MCDA and \cite{Brans_book,Brans84} for a survey on PROMETHEE methods) and widely used in order to deal with several types of real world problems \cite{behzadian2010promethee}. The basic assumption of the PROMETHEE methods, the most known of which are PROMETHEE I and II, is that the set of criteria is mutually preferentially independent \cite{wakker1989additive}. In several decision making problems, this assumption is in general violated since criteria are frequently interacting with each other, which means that it is possible to observe a certain degree of redundancy or synergy among them. For example, when evaluating sport cars, the criteria maximum speed and acceleration may be considered redundant because, in general, speed cars have also a good acceleration. Therefore, even if these two criteria can be very important for a DM liking sport cars, their comprehensive importance is smaller than the sum of importance of the two criteria considered separately. From the other side, criteria maximum speed and price lead to a synergy effect because a speed car having also a low price is very well appreciated. For such a reason, the comprehensive importance of these two criteria should be greater than the sum of importance of the two criteria considered separately. In these two cases the aggregation of the preferences is done by using non-additive integrals the most know of which  are the Choquet integral \cite{choquet1953theory} and the Sugeno integral \cite{sugeno1974theory} (see \cite{Grabisch1996,Grabisch2008,Grabisch_book_greco} for the application of non-additive integrals in MCDA). \\
In other cases, the importance of a criterion expressing preference in favour of an alternative can be reduced by the presence of another criterion opposing to it. For instance, when considering again the car evaluation problem, the weight of criterion speed in favor of an alternative $a$ can be reduced by the presence of criterion aesthetic opposing to it. In these cases, the generalization of Choquet and Sugeno integrals, that are the bipolar Choquet integral and the bipolar Sugeno integral can be used to aggregate the preferences \cite{GL1,GL2,GrecoMatarazzoSlowinski02,GR12}. 

In \cite{Corrente2012a} and \cite{Corrente_annals}, the authors presented the bipolar PROMETHEE methods extending the PROMETHEE I and II methods to the case in which criteria are not independent, proposing to use the bipolar Choquet integral to aggregate the preferences.

In order to apply these methods,  the Decision Maker (DM) has to provide inter-criteria and intra-criterion preference information on the parameters involved. This type of preference information can be given directly or indirectly. Direct preference information consists of asking the DM to provide precise values of the parameters, while indirect preference information consists of asking the DM to provide some comparisons regarding reference alternatives from which the values of the parameters should be elicited. The indirect preference information demands less cognitive effort from the DM and, for this reason, it is widely used in outranking methods (see for example \cite{dias2006inferring,mousseau1998inferring,mousseauuser,GKMS2010,Promethee2010}).  

Recommendations provided by the PROMETHEE methods are dependent on the values given to the considered parameters; in fact, generally, different sets of parameters lead to different comparisons between the alternatives. For this reason, many methods have been proposed in literature to elicit sets of parameters compatible with some preference information provided by the DM or to get values of the parameters giving to an alternative the best rank. In the following we cite only some of these contributions. \\
Considering an MCDA additive method, and getting consequently a ranking of the considered alternatives, Mareschal \cite{mareschal} makes a sensitivity analysis to obtain stability intervals for the criteria weights such that the method provides the same ranking whichever the choice of the weights inside these intervals. Sun and Han \cite{sun} solve a linear programming problem to find the most discriminant set of weights compatible with the preference information provided by the DM. Wolters and Mareschal \cite{wolters} propose a method to study how the ranking of the alternatives is sensitive to the changes on the weights and on the alternatives' evaluations, and which modifications on the importance weights are necessary to make an alternative the best. Eppe and De Smet \cite{eppe_IPMU} study the weights that best represent the preferences of the DM and the type of information structure giving the ranking closest to the ranking provided by the DM Eppe et al \cite{eppe_de_Smet}, use a meta-heuristic in order to induce parameters for the PROMETHEE II model starting from some preference information provided by the DM.

If the DM decides to give indirect preferences, Robust Ordinal Regression (ROR) \cite{greco2008ordinal,figueira2009building,greco2010robust,GKMS2010} can be applied to explore the whole space of parameters compatible with the preference information provided by the DM. ROR is a family of MCDA methods taking into account simultaneously all the sets of parameters compatible with some preference information provided by the DM considering two preference relations, one possible and one necessary. The necessary preference relation holds between two alternatives $a,b$ if $a$ is at least as good as $b$ for all sets of parameters compatible with the information provided by the DM, while the possible preference relation between two alternatives holds if $a$ is at least as good as $b$ for at least one set of parameters compatible with the preference information provided by the DM. ROR for preferences with interacting criteria modelled through the Choquet integral has been considered in \cite{angilella2010non}. ROR has been applied also to the classical PROMETHEE methods in \cite{Promethee2010} and to the bipolar PROMETHEE methods in \cite{Corrente_annals}.

Similarly to the ROR, the Stochastic Multiobjective Acceptability Analysis (SMAA) \cite{Lahdelma,Lahdelma_S2} explores the whole set of parameters compatible with some preference information provided by the DM. SMAA is a family of MCDA methods taking into account imprecision or lack of the data considering probability distributions on the space of criteria weights and on the space of alternatives' evaluations. SMAA has been applied to preferences with interacting criteria modelled through the Choquet integral in \cite{angilella_IPMU}.\\
In this paper we propose an integrated approach between SMAA and the classical and the bipolar PROMETHEE methods giving arise to the SMAA-PROMETHEE methods. From one side, integrating SMAA and PROMETHEE I methods permits to compute the frequency of the preference, indifference or incomparability between two alternatives From the other side, integrating SMAA and PROMETHEE II methods permits to study the different positions an alternative can get when varying the parameters of the model.

The paper is organized as follows. In the next Section we shall describe the classical PROMETHEE I and II methods while in Section \ref{bipol_PROM} we recall the basic concepts of the bipolar PROMETHEE methods; Section \ref{smaa_methods} contains the description of SMAA methods when the preference model is an additive function; in Section \ref{SMAA_PROM} we shall present the SMAA-PROMETHEE methods and we shall describe its link with the ROR; Section \ref{example} contains a didactic example while some conclusions and future directions of research are presented in the last section.

\section{The classical PROMETHEE methods}\label{class_PROM}
Let us consider a set of actions or alternatives $A$ evaluated with respect to a set of criteria $G=\left\{g_{1},\ldots,g_{n}\right\}$, where $g_{j}:A\rightarrow\rio$, $j\in{\cal J}=\left\{1,\ldots,n\right\}$ and $|A|=m$. PROMETHEE \cite{Brans_book,BransVi85} 
is a well-known family of MCDA methods that aggregate preference information of a DM through an outranking relation. Considering for each criterion $g_{j}$ a weight $w_{j}$ (representing the importance of criterion $g_{j}$ within the family of criteria $G$), an indifference threshold $q_{j}$ (being the largest difference $d_{j}(a,b)=g_{j}(a)-g_{j}(b)$ compatible with the indifference between alternatives $a$ and $b$ on criterion $g_j$), and a preference threshold $p_{j}$ (being the minimum difference $d_{j}(a,b)$ compatible with the preference of $a$ over $b$ on criterion $g_j$), PROMETHEE methods build a non-decreasing function of $d_{j}(a,b)$, whose formulation (see \cite{Brans_book} for alternative formulations) can be stated as follows

\begin{equation}\label{p_j}
P_{j}(a,b)=
\left\{
\begin{array}{lll}
0 & \mbox{if} & d_{j}(a,b)\leq q_{j}, \\
\frac{d_{j}(a,b)-q_j}{p_{j}-q_{j}} & \mbox{if} & q_{j}<d_{j}(a,b)<p_{j}, \\
1 & \mbox{if} & d_{j}(a,b)\geq p_{j}.
\end{array}
\right.
\end{equation}

\noindent $P_{j}(a,b)$ represents the degree of preference of $a$ over $b$ on criterion $g_{j}$. \\
For each ordered pair of alternatives $(a,b)\in A\times A,$ PROMETHEE methods compute  
$$
\pi(a,b)=\sum_{j\in{\cal J}}w_{j}P_{j}(a,b)
$$
representing how much alternative $a$ is preferred to alternative $b$ taking into account the whole set of criteria. $\pi(a,b)$ can assume values between $0$ and $1$ and obviously the greater the value of $\pi(a,b)$, the greater the preference of $a$ over $b$. \\
In order to compare an alternative $a$ with all the other alternatives of the set $A$, PROMETHEE methods consider the negative and the positive flow of $a$ defined in the following way:
$$
\phi^{-}(a)=\frac{1}{m-1}\sum_{b\in A\setminus\left\{a\right\}}\pi(b,a) \;\;\;\;\;\;\mbox{and}\;\;\;\;\;\;  \phi^{+}(a)=\frac{1}{m-1}\sum_{b\in A\setminus\left\{a\right\}}\pi(a,b).
$$
\noindent $\phi^{-}(a)$ represents how much the alternatives of $A\setminus\left\{a\right\}$ are preferred to $a$; the lower $\phi^{-}(a)$ the better alternative $a$ is; $\phi^{+}(a)$ represents how much $a$ is preferred to the alternatives of $A\setminus\left\{a\right\}$; the greater $\phi^{+}(a)$, the better $a$ is. PROMETHEE II computes also the net flow $\phi(a)=\phi^{+}(a)-\phi^{-}(a)$ for each alternative $a\in A$.  Taking into account the negative and the positive flows, PROMETHEE I builds a preference (${\cal P}^{I}$), an indifference (${\cal I}^{I}$) and an incomparability (${\cal R}^{I}$) relation on the set of alternatives $A$, while considering the net flow, PROMETHEE II builds a preference (${\cal P}^{II}$) and an indifference (${\cal I}^{II}$) relation on the set of alternatives $A$: \\

$
\left\{
\begin{array}{lll}
a{\cal P}^{I}b  & \mbox{iff} & \Phi^{+}(a)\geq\Phi^{+}(b), \;\Phi^{-}(a)\leq\Phi^{-}(b) \;\;\mbox{and at least one of the two inequalities is strict}, \\[1mm]
a{\cal I}^{I}b  & \mbox{iff} & \Phi^{+}(a)=\Phi^{+}(b) \;\;\mbox{and}\;\; \Phi^{-}(a)=\Phi^{-}(b),\\[1mm]
a{\cal R}^{I}b  & \multicolumn{2}{l}{\mbox{otherwise.}} \\[3mm]
\end{array}
\right.
$

\vspace{0,5cm}
$
\left\{
\begin{array}{lll}
a{\cal P}^{II}b & \mbox{iff} & \Phi(a)>\Phi(b), \\[1mm]
a{\cal I}^{II}b & \mbox{iff} & \Phi(a)=\Phi(b).
\end{array}
\right.
$

\vspace{0,5cm}
\noindent PROMETHEE I and PROMETHEE II provide, respectively, a partial and a complete order of the alternatives in $A$.

\section{The bipolar PROMETHEE methods}\label{bipol_PROM}
In order to take into account the bipolar preferences, the bipolar PROMETHEE methods \cite{Corrente_annals,Corrente2012a} build for each criterion $j\in{\cal J}$, and for each pair of alternatives $(a,b)\in A\times A$, the bipolar preference function

\begin{equation}\label{p_j^B}
P_{j}^{B}(a,b)=
\left\{
\begin{array}{lll}
P_{j}(a,b) & \mbox{if} & P_{j}(a,b)>0 \\
-P_{j}(b,a) & \mbox{if} & P_{j}(a,b)=0.
\end{array}
\right.
\end{equation}

\noindent It is representative of the preference of $a$ over $b$ on criterion $g_j$ and it is straightforward proving that $P_{j}(a,b)=-P_{j}(b,a)$ for all $j\in{\cal J}$ and for all $a,b\in A$.\\
In the bipolar PROMETHEE methods, the aggregation of the preference functions $P_{j}^{B}(a,b)$ is done by using the bipolar Choquet integral. \\
The bipolar Choquet integral is based on a bicapacity \cite{GL1,GL2}, being a function $\hat\mu$ defined on $P({\cal J})=\left\{(C,D): C,D\subseteq {\cal J} \mbox{ and } C \cap D=\emptyset \right\}$ and having value on $\left[-1,1\right]$ such that:
\begin{itemize}
\item $\hat\mu(\emptyset, {\cal J})=-1, \hat\mu({\cal J},\emptyset)=1,$ $\hat\mu(\emptyset, \emptyset)=0$ (boundary constraints),
\item $\hat\mu(C,D)\leq\hat\mu(E,F)$ for all $(C,D),(E,F)\in P({\cal J})$, such that $C\subseteq E$ and $D\supseteq F$ (monotonicity constraints).
\end{itemize}

\noindent According to \cite{GF2003}, in the bipolar PROMETHEE methods we consider a decomposed bicapacity $\hat\mu$ \cite{GrecoMatarazzoSlowinski02} such that $\hat\mu(C,D)=\mu^+(C,D)-\mu^-(C,D)$ for all $(C,D)\in P({\cal J})$ where $\mu^+,\mu^-:P({\cal J})\rightarrow\left[-1,1\right]$ satisfy the following conditions:
\begin{itemize}
\item $\mu^{+}({\cal J},\emptyset)=1$ and $\mu^{+}(\emptyset,B)=0, \;\forall B\subseteq{\cal J}$, 
\item $\mu^{-}(\emptyset,{\cal J})=1$ and $\mu^{-}(B,\emptyset)=0, \;\forall B\subseteq{\cal J},$
\item $\mu^+(C,D) \leq \mu^+(E,F)$ for all $(C,D), \; (E, F) \; \in \; P({\cal{J}}): C \subseteq E, \; D \supseteq F$, 
\item $\mu^-(C,D) \leq \mu^-(E,F)$ for all $(C,D), \; (E, F) \; \in \; P({\cal{J}}): C \supseteq E, \; D \subseteq F$.
\end{itemize}

The interpretation of the functions $\mu^{+}$ and $\mu^{-}$ is the following. Given the pair $(a,b)\in A\times A$, let us consider $(C,D)\in P({\cal J})$ where $C$ is the set of criteria expressing a preference of $a$ over $b$ and $D$ the set of criteria expressing a preference of $b$ over $a$. In this situation, $\mu^{+}(C,D)$ represents the importance of criteria from $C$ when criteria from $D$ are opposing them, and $\mu^{-}(C,D)$ represents the importance of criteria from $D$ opposing $C$. Consequently, $\hat\mu(C,D)$ represents the balance of the importance of $C$ supporting $a$ and $D$ supporting $b$.

After reordering in a non-decreasing way the preference functions for each criterion $j\in{\cal J}$, ($\vert P_{(1)}^{B}(a,b)  \vert \leq \vert P_{(2)}^{B}(a,b) \vert \leq \ldots \leq \vert P_{(j)}^{B}(a,b) \vert \leq \ldots \leq \vert P_{(n)}^{B}(a,b) \vert$), on the basis of the decomposed bicapacity $\hat\mu=\mu^+-\mu^-$, we can give the following definitions of the bipolar comprehensive, positive and negative preference:

\begin{equation}\label{bipolarcomprehensive}
 Ch^{B}(P^{B}(a, b), \hat{\mu}) = \sum_{j \in {\cal{J}}^{>}} \vert P_{(j)}^{B}(a, b) \vert
 \Big[ \hat{\mu}(C_{(j)}, D_{(j)})  - \hat{\mu}(C_{(j+1)}, D_{(j+1)}) \Big],
\end{equation}

\begin{equation}\label{bipolarpositive}
 Ch^{B+}(P^{B}(a, b), {\mu}^{+}) = \sum_{j \in {\cal{J}}^{>}} \vert P_{(j)}^{B}(a, b) \vert
 \Big[ {\mu}^{+}(C_{(j)}, D_{(j)})  - {\mu}^{+}(C_{(j+1)}, D_{(j+1)}) \Big],
\end{equation}

\begin{equation}\label{bipolarnegative}
 Ch^{B-}(P^{B}(a, b), {\mu}^{-}) = \sum_{j \in {\cal{J}}^{>}} \vert P_{(j)}^{B}(a, b) \vert
 \Big[ {\mu}^{-}(C_{(j)}, D_{(j)})  - {\mu}^{-}(C_{(j+1)}, D_{(j+1)}) \Big].
\end{equation}

\noindent where $P^B(a, b) = \Big[P_{j}^{B}(a,b), \; j \in{\cal{J}}\Big]$, ${\cal{J}}^{>} = \{ j \in {\cal{J}} \; : \; \vert P_{(j)}^{B}(a,b) \vert > 0\}$, $C_{(j)} = \{ i \in {\cal{J}}^{>} \; : \; P_{i}^{B}(a,b) \geq \vert P_{(j)}^{B}(a, b) \vert\}$, $D_{(j)} = \{ i \in {\cal{J}}^{>} \; : \; - P_{i}^{B}(a,b) \geq \vert P_{(j)}^{B}(a, b) \vert\}$ and $C_{(n+1)}=D_{(n+1)}=\emptyset$.\\

\noindent $Ch^B(P^B(a,b), \hat{\mu})$ gives the comprehensive
preference of $a$ over $b$ and it is equivalent to $\pi(a,b) -
\pi(b,a)$ in the classical PROMETHEE method. It is, therefore, reasonable expecting that $Ch^B(P^B(a,b), \hat{\mu})=-Ch^B(P^B(b,a), \hat{\mu})$ for all $a,b\in A$.  \\
$Ch^{B+}(P^B(a,b), \mu^+)$ and $Ch^{B-}(P^B(a,b), \mu^-)$ give, instead, how much $a$ outranks $b$ (considering the reasons in favor of $a$) and how much $a$ is outranked by $b$ (considering the reasons against $a$). Also in this case, it is reasonable expecting that $Ch^{B+}(P^B(a,b), \mu^+)=Ch^{B-}(P^B(b,a), \mu^-)$ for all $a,b\in A$. \\
From the definitions above, it is straightforward proving that, for all $a,b\in A,$ 

\begin{equation}\label{bipolarpref}
Ch^{B}(P^{B}(a, b), \hat{\mu})=Ch^{B+}(P^{B}(a, b), {\mu}^{+})-Ch^{B-}(P^{B}(a, b), {\mu}^{-})
\end{equation}

Similarly to the classical PROMETHEE methods, also in the bipolar PROMETHEE methods one can define the bipolar positive, negative and net flows for all $a\in A$:

\begin{equation}\label{pos_flow}
 {\phi}^{B+}(a) = \frac{1}{m-1} \sum_{b \in A\setminus\left\{a\right\}}Ch^{B+}(P^{B}(a,b), {\mu}^{+})
\end{equation}

\begin{equation}\label{neg_flow}
 {\phi}^{B-}(a) = \frac{1}{m-1} \sum_{b \in A\setminus\left\{a\right\}}Ch^{B-}(P^{B}(a,b), {\mu}^{-})
\end{equation}

\begin{equation}\label{net_flow}
 {\phi}^{B}(a) = \frac{1}{m-1} \sum_{b \in A\setminus\left\{a\right\}}Ch^{B}(P^{B}(a,b), \hat{\mu})
\end{equation}

\noindent and by equation (\ref{bipolarpref}), it follows that ${\phi}^{B}(a)={\phi}^{B+}(a)-{\phi}^{B-}(a)$ for all $a\in A.$\\
Equivalently to PROMETHEE I method, the bipolar PROMETHEE I method builds a preference (${\cal P}_{B}^{I}$), an indifference (${\cal I}_{B}^{I}$), and an incomparability (${\cal R}_{B}^{I}$) relation, while equivalently to the classical PROMETHEE II method, the bipolar PROMETHEE II method builds a preference (${\cal P}_{B}^{II}$) and an indifference (${\cal I}_{B}^{II}$) relation on the set of alternatives $A$:\\

$
\left\{
\begin{array}{lll}
a{\cal P}_B^{II}b  & \mbox{iff} & \Phi^{B+}(a)\geq\Phi^{B+}(b), \;\Phi^{B-}(a)\leq\Phi^{B-}(b) \;\;\mbox{and at least one of the two inequalities is strict}, \\[1mm]
a{\cal I}_B^{II}b  & \mbox{iff} & \Phi^{B+}(a)=\Phi^{B+}(b) \;\;\mbox{and}\;\; \Phi^{B-}(a)=\Phi^{B-}(b),\\[1mm]
a{\cal R}_B^{II}b  & \multicolumn{2}{l}{\mbox{otherwise.}} \\[3mm]
\end{array}
\right.
$

\vspace{0,5cm}
$
\left\{
\begin{array}{lll}
a{\cal P}_B^{II}b & \mbox{iff} & \Phi^B(a)>\Phi^B(b), \\[1mm]
a{\cal I}_B^{II}b & \mbox{iff} & \Phi^B(a)=\Phi^B(b).
\end{array}
\right.
$

\vspace{0,5cm}

Because the use of the bipolar Choquet integral, and therefore of the bipolar PROMETHEE methods, is based on a
bicapacity which assigns numerical values to each element of $P({\cal{J}})$, in \cite{Corrente_annals} we have considered only the $2$-additive bicapacities \cite{GL1,fujimoto2004new}, being a particular class of bicapacities. We give the following decomposition of the functions $\mu^{+}$ and $\mu^{-}$ previously defined:

\begin{itemize}
\item ${\displaystyle \mu^{+}(C, D) = \sum_{j \in C}a^{+}(\{ j\}, \emptyset) +  \sum_{\{j,k \} \subseteq C}a^{+}(\{j, k\}, \emptyset) +
  \sum_{j \in C, \; k \in D}a^{+}(\{ j\}, \{ k \}) }$
\item ${\displaystyle \mu^{-}(C, D) = \sum_{j \in D}a^{-}(\emptyset, \{ j\}) +  \sum_{\{j,k \} \subseteq D}a^{-}(\emptyset, \{j, k\}) +
  \sum_{j \in C, \; k \in D}a^{-}(\{ j\}, \{ k \}) }$
\end{itemize}

\noindent where:
\begin{itemize}
\item $a^+(\{ j\}, \emptyset)$ represents the power of criterion $g_j$ by itself expressing positive preferences; this value is always non negative;
\item $a^+(\{j, k\}, \emptyset)$ represents the interaction between $g_j$ and $g_k$ when they are in favor of the preference of $a$ over $b$; this value is zero if there is no interaction; it is positive (negative) when there is a synergy (redundancy) effect when putting together $g_j$ and $g_k$;
\item  $a^+(\{ j\}, \{ k \})$, represents the power of criterion $g_k$ against criterion $g_j$, when criterion
$g_j$ is in favor of $a$ over $b$ and $g_k$ is against to the preference of $a$ over $b$; this leads always to a reduction or no effect on the
value of $\mu^+$ since this value is always non-positive.
\end{itemize}
An analogous interpretation can be applied to the values $a^-(\emptyset, \{
j\})$, $a^-(\emptyset, \{j, k\})$, and $a^-(\{ j \}, \{ k \})$.

Let us remember that positive and negative interactions, and power of one criterion against another correspond to the mutual strengthening effect, the mutual weakening effect and the antagonostic effect considered for the ELECTRE methods in \cite{FGR}.

From $\mu^+$ and $\mu^-$ we get the following \textit{2-additive decomposable bicapacity} $\hat{\mu}$

\begin{equation}
\hat{\mu}(C,D) = \sum_{j \in C}a^{+}_{j} -  \sum_{j \in D}a^{-}_{j} + \sum_{\{j, k \} \subseteq C}a^{+}_{jk} - \sum_{\{j, k \} \subseteq D}a^{-}_{jk} + \sum_{j \in C, \; k \in D}a^{+}_{j \vert k} - \sum_{j \in C, \; k \in D}a^{-}_{j \vert k}
\end{equation}

\noindent where $a_{j}^+$, $a_{jk}^+$, $a_{j \vert k}^+$ denote $a^+(\{ j\},\emptyset)$, $a^+(\{j, k\}, \emptyset)$, $a^+(\{ j\}, \{
k \})$, and $a_{j}^-$, $a_{jk}^-$, $a_{j \vert k}^-$ denote $a^-(\emptyset, \{ j\})$, $a^-(\emptyset, \{j,k\})$, $a^-(\{ j \}, \{ k \})$, respectively.\\

On the basis of a 2-additive decomposable bicapacity, the monotonicity and the boundary conditions previously considered can be rewritten as follows (see \cite{Corrente_annals}):

\begin{itemize}
\item $\displaystyle a^{+}_{j} + \sum_{k \in C}a^{+}_{jk} + \sum_{k \in D}a^{+}_{j \vert
  k} \geq 0, \; \; \forall \; j \in {\cal{J}}, \; \forall  (C \cup \{ j \}, D) \in
  P({\cal{J}})$,
\item $\displaystyle  \sum_{h \in C}a^{+}_{h \vert
  j} \leq 0, \; \; \forall \; j \in {\cal{J}}, \; \forall  (C, D \cup \{ j \}) \in
  P({\cal{J}}),$
\item $\displaystyle a^{-}_{j} + \sum_{k \in D}a^{-}_{jk} + \sum_{h \in C}a^{-}_{h \vert
  j} \geq 0, \; \; \forall \; j \in {\cal{J}}, \; \forall  (C , D\cup \{ j \}) \in
  P({\cal{J}}),$
\item $\displaystyle  \sum_{k \in D}a^{-}_{j \vert
  k} \leq 0, \; \; \forall \; j \in {\cal{J}}, \; \forall  (C\cup \{ j \}, D ) \in
  P({\cal{J}})$,
  \item ${\displaystyle \sum_{j \in {\cal{J}}}a_{j}^{+} +
\sum_{\{j, k \} \subseteq {\cal{J}}}a_{jk}^{+} = 1}$,
  \item ${\displaystyle \sum_{j \in {\cal{J}}}a_{j}^{-} +
\sum_{ \{j, k \} \subseteq {\cal{J}}}a_{jk}^{-} = 1}$.
\end{itemize}

\noindent Considering a 2-additive decomposable bicapacity, in \cite{Corrente_annals} we have given the following theorems defining the comprehensive bipolar positive flow and the comprehensive bipolar negative flow; before we have provided conditions so that $Ch^B(P^B(a,b), \hat{\mu})=-Ch^B(P^B(b,a), \hat{\mu})$ and $Ch^{B+}(P^B(a,b), \mu^+)=Ch^{B-}(P^B(a,b), \mu^-)$ for all $a,b\in A$.

\begin{theorem}(see \cite{Corrente_annals})\label{Bi-polar_Choq} 
{\it Given a 2-additive decomposable bicapacity $\hat{\mu}$, then for all $x \in {\rio }^n$
\begin{enumerate}

\item ${\displaystyle Ch^{B+}(x, {\mu}^{+}) = \sum_{j \in {\cal{J}}, x_j> 0}a^{+}_{j}x_j + \sum_{j, k \in {\cal{J}}, j \neq k, x_j, x_k > 0}a^{+}_{jk}\min\{x_j,
 x_k\}+ \sum_{j, k \in {\cal{J}}, j\neq k, x_j > 0, x_k < 0}a^{+}_{j \vert k}\min\{x_j, - x_k\} }
$

\item ${\displaystyle Ch^{B-}(x, {\mu}^{-}) = -\sum_{j \in {\cal{J}}, x_j< 0}a^{-}_{j}x_j - \sum_{j, k \in {\cal{J}}, j \neq k, x_j, x_k < 0}a^{-}_{jk}\max\{x_j,
 x_k\}- \sum_{j, k \in {\cal{J}}, j\neq k, x_j > 0, x_k < 0}a^{-}_{j \vert k}\max\{-x_j, x_k\} }$
\end{enumerate}
}
\end{theorem}

\begin{theorem}(see \cite{Corrente_annals})\label{lem_SYM}
Given a 2-additive decomposable bicapacity $\hat\mu$, $Ch^{B}(P^{B}(a,b),\hat\mu)=-Ch^{B}(P^{B}(b,a),\hat\mu)$ and $Ch^{B+}(P^B(a,b), \mu^+)=Ch^{B-}(P^B(b,a), \mu^-)$ for all $a,b\in A$ iff 
\begin{enumerate}
\item for each $j \in {\cal{J}}$, $a^{+}_{j} = a^{-}_{j}$,
\item for each $\{j,k\} \subseteq {\cal{J}}$, $a^{+}_{jk} = a^{-}_{jk}$,
\item for each $ j, k \in {\cal{J}}$, $j \neq k$, $a^{+}_{j \vert k} = a^{-}_{k \vert j}$. 
\end{enumerate}
\end{theorem}

\section{SMAA}\label{smaa_methods}
Chosen the decision model and assigned the values of its parameters, it is straightforward obtaining some recommendation regarding the problem at hand. In general, the DM is not able to provide directly these parameters or she is not worth of their meaning; besides, in real world decision making problems, data regarding the alternatives' evaluations can be missing or imprecise. 

SMAA is a family of MCDA methods taking into account imprecision or lack of data in the problem at hand (see \cite{tervonen_figueira} for a comprehensive survey on the use of SMAA methods in MCDA) considering a probability distribution $f_{W}$ over the space of all compatible weights $W$ and a probability distribution $f_{\chi}$ over the space $\chi\subseteq\rea^{m\times n}$ of the alternatives' evaluations. \\
Without any preference information provided by the DM, the weight space is defined as follows: 

$$\left\{w\in\rea^{n}: w_j\geq 0, \forall j\in{\cal J},\; \mbox{and} \; \sum_{j\in{\cal J}}w_j=1\right\}.$$ 

\noindent If the DM is able to provide some preference information, then the space $W$ is restricted by the constraints translating these preference information. Consequently, in the following, with $W$ we shall denote the set of weights compatible with the preference information provided by the DM. Now let us describe some SMAA concepts considering the additive value function, $\displaystyle u(x_i,w)=\sum_{j\in{\cal J}}w_j g_j(x_i)$ with $x_i\in A$ and $w\in W$, as preference model.

Because $u(x_i,w)$ provides a complete ranking of the alternatives, for each $\xi$ in $\chi$ and $w$ in $W$, SMAA computes the position reached by alternative $x_i\in A$ as 

$$rank(i,\xi,w)=1+\sum_{k\neq i}\rho(u(\xi_k,w)>u(\xi_i,w)),$$

\noindent where $\rho(true)=1$ and $\rho(false)=0.$ Besides, for each $\xi\in\chi$, SMAA computes the favourable rank weights of alternative $x_i\in A$

$$W_{i}^{r}(\xi)=\left\{w\in W:rank(i,\xi,w)=r\right\}$$
 
\noindent being the set of possible weights giving to alternative $x_i$ the position $r=1,\ldots,m$ in the final rank. \\
On the basis of the favourable rank weights, SMAA computes the rank acceptability index and the central weight vector.

\begin{itemize}
\item The rank acceptability index

$$b_{i}^{r}=\int_{\xi\in \chi}f_{X}(\xi)\int_{w\in W_{i}^{r}(\xi)}f_{W}(w)\;dw\;d\xi$$

\noindent describes the share of parameters giving to alternative $x_i$ the position $r$ in the obtained final ranking; in particular, $b_{i}^1$ measures the variety of parameters making alternative $x_i$ the most preferred. Obviously, the best alternatives are those having rank acceptability index greater than zero for the first positions and rank acceptability index close to zero for the lower positions. 
\item The central weight vector

$$w_{i}^{c}=\frac{1}{b_{i}^{1}}\int_{\xi \in \chi}f_{X}(\xi)\int_{w\in W_{i}^{1}(\xi)}f_{W}(w)w\;dw\;d\xi$$

\noindent describes the preferences of a typical DM that makes alternative $x_i$ the most preferred.

%

\end{itemize} 

\section{SMAA-PROMETHEE}\label{SMAA_PROM}
Even if SMAA has already been applied to outranking methods and in particular to ELECTRE methods \cite{hokkanen,tervonen_electre_tri}, to the best of our knowledge, no attempt of applying the SMAA methods to PROMETHEE has been done in literature before. In \cite{Corrente_annals}, we pointed out that the bipolar PROMETHEE methods are reduced to the classical PROMETHEE methods if $a_{jk}=a_{j\vert k}^{+}=a_{j\vert k}^{-}=0$ for all $j,k\in{\cal J}, j\neq k$ and in this case the parameters $a_j$ in the bipolar PROMETHEE methods are the weights $w_j$ of the classical PROMETHEE methods. For this reason, without loss of generality, in this section we shall apply SMAA only to the bipolar PROMETHEE methods. 

We shall suppose that parameters $a_{j},$ $a_{j\vert k}^+,$ and $a_{j\vert k}^-,$ are unknown for all $j,k\in{\cal J}$ and, at the same time, we shall assume that indifference and the preference thresholds as well as the alternatives' evaluations are already known.

In the bipolar PROMETHEE methods, we consider the following information
provided by the DM and their representation in terms of linear
constraints:

\begin{enumerate}
\item {\it Comparing pairs of actions locally or globally}. The constraints represent some pairwise comparisons on a set of
training actions. Given two actions $a$ and $b$, the DM may
prefer $a$ to $b$, $b$ to $a$ or be indifferent to both:
    \begin{enumerate}
    \item the linear constraint associated with $a{\cal P}b$ ($a$ is locally preferred to $b$) is: 
    $$Ch^{B}(P^{B}(a, b), \hat{\mu}) > 0;$$
    \item the linear constraints associated with $a{\cal P}_B^{I}b$ ($a$ is preferred to $b$ with respect to the bipolar PROMETHEE I method) are:     
    $$    
    \left.
    \begin{array}{l}
    \Phi^{B+}(a)\geq\Phi^{B+}(b),\\
    \Phi^{B-}(a)\leq\Phi^{B-}(b),\\
    \Phi^{B+}(a)-\Phi^{B-}(a)>\Phi^{B+}(b)-\Phi^{B-}(b),\\
    \end{array}
    \right\} 
    $$
    \item the linear constraint associated with $a{\cal P}_B^{II}b$ ($a$ is preferred to $b$ with respect to the bipolar PROMETHEE II method) is: 
    $$    
    \Phi^{B}(a)>\Phi^{B}(b)
    $$    
    \item the linear constraint associated with $a{\cal I}b$ ($a$ is locally indifferent to $b$) is:
    $$
    Ch^{B}(P^{B}(a, b), \hat{\mu}) = 0
    $$
    \item the linear constraints associated with $a{\cal I}_B^{I}b$ ($a$ is indifferent to $b$ with respect to the bipolar PROMETHEE I method) are:    
    $$    
    \left.
    \begin{array}{l}
    \Phi^{B+}(a)=\Phi^{B+}(b),\\
    \Phi^{B-}(a)=\Phi^{B-}(b),\\
    \end{array}
    \right\}
    $$
    \item the linear constraint associated with $a{\cal I}_B^{II}b$ ($a$ is indifferent to $b$ with respect to the bipolar PROMETHEE II method) is:
    $$
    \Phi^{B}(a)=\Phi^{B}(b)
    $$
    \end{enumerate}
\item {\it Comparison of the intensity of preferences between pairs of
actions}. The constraints represent some pairwise comparisons between pairs of alternatives on a set of
training actions. Given four actions $a$, $b$, $c$ and $d$:
\begin{enumerate}
\item the linear constraints associated with $(a,b) {\cal{P}} (c, d)$ (the
local preference of $a$ over $b$ is larger than the
local preference of $c$ over $d$) is:
\[
Ch^{B}(P^{B}(a, b), \hat{\mu}) > Ch^{B}(P^{B}(c, d), \hat{\mu})
\]
\item the linear constraints associated with $(a,b) {\cal{I}} (c, d)$ (the
local preference of $a$ over $b$ is the same of
local preference of $c$ over $d$) is:
\[
Ch^{B}(P^{B}(a, b), \hat{\mu}) = Ch^{B}(P^{B}(c, d), \hat{\mu}) 
\]
\end{enumerate}
\item {\it Importance of criteria}.
A partial ranking over the set of criteria $\cal{J}$ may be
provided by the DM: 
    \begin{enumerate}
    \item criterion $g_j$ is more important than criterion $g_k$, which leads
    to the constraint $a_{j} > a_{k}$;
    \item criterion $g_j$ is equally important to criterion $g_k$, which leads
    to the constraint $a_{j} = a_{k}$.
    \end{enumerate}
\item {\it The sign of interactions}. The DM may be able, for certain cases,
to provide the sign of some interactions. For example, if there is
a synergy effect when criterion $g_j$ interacts with criterion
$g_k$, the following constraint
should be added to the model: $a_{jk} > 0$. 
\item {\it Interaction between pairs of criteria}. The DM can provide
some information about interaction between criteria:
     \begin{enumerate}
    \item[a)] if the DM feels that interaction between $g_j$ and $g_k$ is greater
    than the interaction between $g_p$ and $g_q$, the constraint should be
    defined as follows: $|a_{jk}| > |a_{pq}|$ where in particular:
    \begin{itemize}
    \item if both couples of criteria are synergic then: $a_{jk} > a_{pq}$,
    \item if both couples of criteria are redundant then: $a_{jk} < a_{pq}$,
    \item if $(j,k)$ is a couple of synergic criteria and $(p,q)$ is a couple of redundant criteria, then: $a_{jk} > -a_{pq}$,
    \item if $(j,k)$ is a couple of redundant criteria and $(p,q)$ is a couple of synergic criteria, then: $-a_{jk} > a_{pq}$.
    \end{itemize}
    \item[b)] if the DM feels that the strength of the interaction between $g_j$ and $g_k$ is the same of the strength of the interaction between
    $g_p$ and $g_q$, the constraint will be the following: $|a_{jk}| = |a_{pq}|$ and in particular:
    \begin{itemize}
    \item if both couples of criteria are synergic or redundant then: $a_{jk} = a_{pq}$,
    \item if one couple of criteria is synergic and the other is redundant then: $a_{jk} = -a_{pq}$,
    \end{itemize}     
    \end{enumerate}
\item {\it The power of the opposing criteria}{\label{interact}}. Concerning the power
of the opposing criteria several situations may occur. For
example:
    \begin{enumerate}
    \item[a)] when the opposing power of $g_k$ is larger than the
    opposing power of $g_h$, with respect to $g_j$, which expresses a positive preference,
    we can define the following constraint:
    $a^{+}_{j \vert k} < a^{+}_{j \vert h}$ (because $a^{+}_{j \vert h}\leq 0$ and $a^{-}_{j \vert h}\leq 0$ for all $j,k$ with $j\neq k$);
    \item[b)] if the opposing power of $g_k$, expressing negative preferences, is larger with $g_j$ rather
    than with $g_h$, the constraint will be $a^{+}_{j \vert k} < a^{+}_{h \vert
    k}$.
    \end{enumerate}
\end{enumerate}

The constraint translating the preference information of the DM along with the symmetry, boundary and monotonicity conditions make the following set of constraints in which strict inequalities have been transformed in weak inequalities by adding an auxiliary variable $\varepsilon$:

\begin{scriptsize}
    \begin{displaymath}
      \left.
      \begin{array}{ll}
           
             Ch^{B}(P^{B}(a, b),\hat\mu) \geq \varepsilon \;\; \mbox{if} \;\; a{\cal{P}}b, & Ch^{B}(P^{B}(a, b),\hat\mu) = 0  \;\; \mbox{if} \;\; a{\cal{I}}b, \\[1mm]          
             \left.
             \begin{array}{l}
             \Phi^{B+}(a)\geq\Phi^{B+}(b),\\
             \Phi^{B-}(a)\leq\Phi^{B-}(b),\\
             \Phi^{B+}(a)-\Phi^{B-}(a)\geq \Phi^{B+}(b)-\Phi^{B-}(b)+\varepsilon\\
             \end{array}
             \right\} \;\; \mbox{if} \;\; a{\cal P}_B^{I}b & 
             \left.
             \begin{array}{l}
             \Phi^{B+}(a)=\Phi^{B+}(b),\\
             \Phi^{B-}(a)=\Phi^{B-}(b)\\           
             \end{array}
             \right\} \;\; \mbox{if} \;\; a{\cal I}_B^{I}b \\[1mm]
           \Phi^{B}(a)\geq\Phi^{B}(b)+\varepsilon \;\; \mbox{if} \;\; a{\cal P}_B^{II}b & \Phi^{B}(a)=\Phi^{B}(b) \;\; \mbox{if} \;\; a{\cal I}_B^{II}b \\[1mm]
             {Ch^{B}(P^{B}(a, b),\hat\mu) \geq Ch^{B}(P^{B}(c, d),\hat\mu) + \varepsilon } \;\; \mbox{if} \;\; (a,b) {\cal{P}} (c,d), & {Ch^{B}(P^{B}(a, b),\hat\mu) = Ch^{B}(P^{B}(c, d),\hat\mu)}  \;\; \mbox{if} \;\;(a,b){\cal{I}}(c,d), \\[1mm]
            a_{j} - a_{k} \geq \varepsilon \;\; \mbox{if} \;\; j{\cal{P}}k,  & a_{j} = a_{k} \;\; \mbox{if} \;\; j{\cal{I}}k,   \\[1mm] 
            |a_{jk}| - |a_{pq}| \geq \varepsilon \;\; \mbox{if} \;\;  \{ j, k\} {\cal{P}}\{ p, q\}, \; \mbox{(see point 5.a) above)}& \\[1mm]
            |a_{jk}| = |a_{pq}|                \;\;\mbox{if} \;\;  \{ j, k\} {\cal{I}}\{ p, q\}, \; \mbox{(see point 5.b) above)} &  \\[1mm]  
            a_{jk} \geq \varepsilon  \;\; \mbox{if there is synergy between criteria $j$ and $k$}, \\[1mm]
            a_{jk} \leq - \varepsilon \;\; \mbox{if there is redundancy between criteria $j$ and $k$}, \\[1mm]
            a_{jk} = 0  \;\;\mbox{if criteria $j$ and $k$ are not interacting}, \\[1mm]
           \mbox{Power of the opposing criteria of the type \ref{interact}:}\\[1mm]
           a^{+}_{j \vert k} - a^{+}_{j \vert p} \geq \varepsilon, & a^{-}_{j \vert k} - a^{-}_{j \vert p} \geq \varepsilon,  \\[1mm]
           a^{+}_{j \vert k} - a^{+}_{p \vert k} \geq \varepsilon, & a^{-}_{j \vert k} - a^{-}_{p \vert k} \geq \varepsilon,  \\[1mm]           
           \mbox{Symmetry conditions (Theorem \ref{lem_SYM}):}\\[1mm] 
           {\displaystyle  a_{j \vert k}^{+} = a_{k \vert j}^{-}, \; \; \forall \; j, k\in {\cal{J}}, j\neq k \;\; \; }\\[1mm]
           \mbox{Boundary and monotonicity conditions:}\\[1mm] 
           {\displaystyle \sum_{j \in {\cal{J}}}a_{j} + \sum_{\{j, k \} \subseteq {\cal{J}}}a_{jk} = 1},\\[1mm]
           {\displaystyle a_{j}\geq 0 \; \; \forall \; j \in {\cal{J}}}, & {\displaystyle a_{j \vert k}^{+}, \; a_{j \vert k}^{-} \; \leq 0 \; \; \forall \; j, k \in {\cal{J}}},\\[1mm]
           {\displaystyle a_{j} + \sum_{k \in C}a_{jk} + \sum_{k \in D}a^{+}_{j \vert k} \geq 0, \; \; \forall \;
           j \in {\cal{J}}, \; \forall  (C \cup \{ j \}, D) \in P({\cal{J}}) }, \\[1mm]
            {\displaystyle a_{j} + \sum_{k \in D}a_{jk} + \sum_{h \in C}a^{-}_{h \vert j} \geq 0, \; \; \forall \;
           j \in {\cal{J}}, \; \forall  (C, D \cup \{ j \}) \in P({\cal{J}}) }. \\           
    \end{array}
    \right\}E^{B}
    \end{displaymath}
\end{scriptsize}

In order to progress from the simplest to the most sophisticated method, one could check firstly if the preferences of the DM are compatible with the classical PROMETHEE methods; if this is not true, then one can pass to the bipolar PROMETHEE methods. 

\begin{note}
The DM could start directly by checking if the bipolar PROMETHEE methods can restore her preference information without looking at the compatibility of the classical PROMETHEE methods; in the same way, she can check if the bipolar PROMETHEE methods can restore her preferences even if the classical PROMETHEE methods are already able to restore them.
\end{note}
In order to check if the preference information provided by the DM can be restored by the classical PROMETHEE methods, one has to solve the following optimization problem:

\begin{equation}\label{first_probl}
\begin{array}{l}
\;\;\mbox{max}\; \varepsilon=\varepsilon_1^{*}\\[2mm]
\left.
\begin{array}{l}
\; E^{B}\\[1mm]
\; a_{jk}=a^{+}_{j\vert k}=a^{-}_{j\vert k}=0, \;\forall j,k\in{\cal J}
\end{array}
\right\}E^{B'}
\end{array}
\end{equation}

\noindent If $E^{B'}$ is feasible and $\varepsilon_1^{*}>0$, then the classical PROMETHEE methods are able to restore the preference information provided by the DM, otherwise, solve the following optimization problem

\begin{equation}\label{second_prob}
\begin{array}{l}
\;\;\mbox{max}\; \varepsilon=\varepsilon_2^{*} \;\;\;s.t.\\[1mm]
\;\; E^{B}\\[1mm]
\end{array}
\end{equation}

\noindent If $E^{B}$ is feasible and $\varepsilon_2^{*}>0$, then the bipolar PROMETHEE methods can restore the preference information provided by the DM, otherwise one can try to discover which constraints cause the incompatibility by using some technique presented in \cite{mousseau2003resolving}.

Since the constraints in $E^{B}$ are linear, we propose to apply the Hit-And-Run method \cite{smith1984,Tervonen2012} to sample sets of parameters compatible with the preference information provided by the DM.\\
The Hit-And-Run sampling begins with the choice of one point inside the polytope delimited by the constraints translating the preference information provided by the DM. At each iteration, it is sampled from the unit hypersphere a random direction that, with the considered position, generate a line. Finally, the sampling of one point inside the segment whose extremes are the intersection of the line with the bounds and the current point is done. \\
Proceeding in this way, we sample a given number of sets of parameters and, applying the PROMETHEE methods considering each time one of these compatible sets of parameters, we store:

\begin{enumerate}
\item preference ${\cal P}^{I}_{B}$, indifference ${\cal I}^{I}_{B}$ and incomparability ${\cal R}^{I}_{B}$ relations obtained by applying the bipolar PROMETHEE I method,
\item preference ${\cal P}^{II}_{B}$ and indifference ${\cal I}^{II}_{B}$ relations obtained by applying the bipolar PROMETHEE II method,
\item ranking of the alternatives by applying the bipolar PROMETHEE II method,
\item set of parameters sampled in the iteration.
\end{enumerate}

\noindent Sampled the fixed number of sets of parameters, we get the following results:

\begin{itemize}
\item Considering the ranking obtained by the bipolar PROMETHEE II method for each parameters' sampling, one can compute the rank acceptability index $b_{i}^r$, the preference matrix $P$ and the indifference matrix $I$. In particular, $b_{i}^r$ describes the share of parameter values giving to alternative $x_i$ rank $r$; the element $p_{ij}$ of the preference matrix $P$ is the frequency of the preference of alternative $x_i$ over $x_j$ while the element of the indifference matrix $I$ in the position $(i,j)$ is the frequency of the indifference between $x_i$ and $x_j$. Besides, one can compute the central weight vector $w_i^{c}$ and the barycenter $w^c$ representing respectively the typical preferences allowing to alternative $x_i$ to get the first rank and the average preferences of the DM;
\begin{note}
We observe that $w_{i}^{c}$ and $w_{c}$ are vectors of  $\rea^{n}$ if we apply SMAA to the classical PROMETHEE methods, while they belong to $\rea^{\frac{3n^2-n}{2}}$ if we apply SMAA to the bipolar PROMETHEE methods ($\frac{3n^2-n}{2}=n+C_{n,2}+D_{n,2}$ where $n$ are the weights of criteria, $C_{n,2}$ is the number of the parameters regarding the interaction between couples of criteria and $D_{n,2}$ is the number of parameters of the type $a^{+}_{j\vert k}$. We did not consider the parameters $a^{-}_{j\vert k}$ because $a^{+}_{j\vert k}=a^{-}_{k\vert j}$).
\end{note}
\item Considering the bipolar PROMETHEE I method, for each couple of alternatives $x_i,x_j\in A$ one can compute the frequency of the preference of $x_i$ over $x_j$, the frequency of the preference of $x_j$ over $x_i$ or still the frequency of the indifference or incomparability between $x_i$ and $x_j$;
\end{itemize}

In \cite{Corrente_annals}, the authors have applied the ROR \cite{greco2010robust} to the bipolar PROMETHEE methods. ROR is a family of MCDM methods tacking into account simultaneously all the sets of parameters compatible with preference information provided by the DM building a necessary and a possible preference relation. Given two alternatives $x_i$ and $x_j$, we say that $x_i$ is necessarily preferred to $x_j$, and we write $x_i\succsim^{N}x_j$, if $x_i$ is at least as good as $x_j$ for all compatible sets of parameters, while $x_{i}$ is possibly preferred to $x_j$, and we have write $x_i\succsim^Px_j$, if $x_i$ is at least as good as $x_j$ for at least one compatible set of parameters. ROR has been already applied to both classical PROMETHEE methods \cite{Promethee2010} and to the bipolar PROMETHEE methods \cite{Corrente_annals}. \\
Because the computation of the necessary and possible preference relations involve to solve two linear programming problems for each pair of alternatives in $A$, we could use the SMAA-PROMETHEE methods to approximate the two preference relations. In fact: 
\begin{itemize}
\item if $x_i$ is necessarily preferred to $x_j$, then the sum of the frequency of the preference of $x_{i}$ over $x_{j}$ and of the frequency of the indifference between $x_{i}$ and $x_{j}$ is 100\%,
\item if the frequency of the preference of $x_i$ over $x_j$ is greater than zero, then $x_i$ is possibly preferred to $x_j$. 
\end{itemize}
The vice versa of these statement are not true because, from one side, also if the frequency of the preference of $x_i$ over $x_j$ is the 100\%, there could exists one non sampled set of compatible parameters for which $x_j$ is preferred to $x_i$; from the other side, also if $x_i$ is possibly preferred to $x_j$, it is possible that for all sampled sets of parameters $x_{j}$ is at least as good as $x_i$ and therefore the frequency of the preference of $x_i$ over $x_j$ is 0\%. Observe that the larger is the sample of the sets of parameters, the better the approximation of the necessary and the possible preference relations obtained through SMAA, such that, in case of an enough large sample of parameters, one can reasonably accept the approximation of SMAA as results of ROR analysis.

\section{Illustrative example}\label{example}
In this section we shall apply the SMAA-PROMETHEE method described in section \ref{SMAA_PROM} to the same example presented in \cite{Corrente_annals}. \\
Suppose that the Dean of a high school has to decide which student deserves to get a scholarship. Eight students arrived to the final selection and in order to make the ``best'' decision, the Dean decides to evaluate them with respect to three different subjects: Mathematics (M), Physics (P) and Literature (L). Their evaluations are given on a [0,20] scale as shown in Table \ref{Evaluations_no_ROR}. The Dean has some preferences regarding two couples of students: $s_7$ and $s_2$ from one side and $s_5$ and $s_6$ from the other side. She states that, she prefers $s_7$ to $s_2$ and $s_5$ to $s_6$. These local preference information are translated by the constraints $Ch^{B}(P^{B}(s_7,s_2),\hat\mu)>0$ and  $Ch^{B}(P^{B}(s_5,s_6),\hat\mu)>0$ respectively. 

\begin{table}[htbp]
\begin{center}
\caption{Evaluations of the students}\label{Evaluations_no_ROR}
\begin{tabular}{cccc}
\hline
\mbox{Students} & \mbox{Mathematics} & \mbox{Physics} & \mbox{Literature}\\
\hline
$s_1$   & 16 & 16 & 16 \\
$s_2$   & 15 & 13 & 18 \\
$s_3$   & 19 & 18 & 14 \\
$s_4$   & 18 & 16 & 15 \\
$s_5$   & 15 & 16 & 17 \\
$s_6$   & 13 & 13 & 19 \\
$s_7$   & 17 & 19 & 15 \\
$s_8$   & 15 & 17 & 16 \\
\hline
\end{tabular}
\end{center}
\end{table}

As explained in the previous section, at first we check if the preferences of the Dean are compatible with the classical PROMETHEE methods. In this case we get that the set of constraints $E^{B'}$ is feasible and $\varepsilon_{1}^{*}>0$. Therefore, the classical PROMETHEE methods can restore the preference information provided by the Dean. \\
Since the Dean is interested to get a general overview about the eight students, we present to her the results in Tables \ref{first_table}, and \ref{second_table} got by applying SMAA to the classical PROMETHEE I and II.

\begin{table}[!h]
\begin{center}
\caption{Results obtained by applying SMAA to the classical PROMETHEE I method}\label{first_table}
\subtable[Preferences in percentage\label{Pref_PROM_I}]{%
\resizebox{0.4\textwidth}{!}{
\begin{tabular}{ccccccccc}
\hline
    &              $\mathbf{s_1}$  &  $\mathbf{s_2}$ &  $\mathbf{s_3}$ &  $\mathbf{s_4}$ &  $\mathbf{s_5}$ &  $\mathbf{s_6}$ &  $\mathbf{s_7}$ &  $\mathbf{s_8}$ \\ 
\hline
 $\mathbf{s_1}$  &        0        &     67,91       &       0         &       18,44     &      58,70     &      61,07    &       0         & 42,96 \\
 
 $\mathbf{s_2}$  &        0        &      0          &     2,09       &         0       &       7,46      &      38,99      &        0        & 11,85 \\
 
 $\mathbf{s_3}$  &     71,25       &    91,21       &       0         &       68,01    &      73,60     &      88,99     &     35,70      & 71,39 \\
 
 $\mathbf{s_4}$  &   63,43 & 81,31 & 0     & 0     & 65,77 & 73,74 & 27,83 & 55,67 \\
 
 $\mathbf{s_5}$  &   31,06 & 62,89 & 2,59  & 26,62 & 0     & 58,09 & 0     & 37,10 \\
 
 $\mathbf{s_6}$  &   0     & 9,52 & 1,75  & 0     & 0     & 0     & 0     & 0 \\
 
 $\mathbf{s_7}$  &   90,90 & 99,19 & 39,23 & 71,08 & 87,60 & 84,83 & 0     & 89,21 \\
 
 $\mathbf{s_8}$  &   47,4  & 62,89 & 1,55 & 38,48 & 62,9 & 60,40 & 0     & 0 \\
\hline
\end{tabular}
}
}
\subtable[Incomparability in percentage\label{Incom_PROM_I}]{%
\resizebox{0.4\textwidth}{!}{
\begin{tabular}{ccccccccc}
\hline
    &  $\mathbf{s_1}$  &  $\mathbf{s_2}$ &  $\mathbf{s_3}$ &  $\mathbf{s_4}$ &  $\mathbf{s_5}$ &  $\mathbf{s_6}$ &  $\mathbf{s_7}$ &  $\mathbf{s_8}$ \\
\hline

 $\mathbf{s_1}$  &        0     & 32,08 & 28,74 & 18,12 & 10,23 & 38,92 & 9,09 & 9,63 \\
 
 $\mathbf{s_2}$  &        & 0     & 6,69  & 18,68 & 29,64 & 51,48 & 0,80 & 25,24 \\
 
 $\mathbf{s_3}$  &    &   & 0     & 31,98 & 23,80 & 9,25 & 25,06 & 27,05 \\
 
 $\mathbf{s_4}$  &    &  &  & 0     & 7,59 & 26,25 & 1,07 & 5,83 \\
 
 $\mathbf{s_5}$  &    &  &  &  & 0     & 41,90 & 12,39 & 0 \\
 
 $\mathbf{s_6}$  &    &  &  &  &  & 0     & 15,16 & 39,59 \\
 
 $\mathbf{s_7}$  &    &  &  &  &  &  & 0     & 10,78 \\
 
 $\mathbf{s_8}$  &    &  &  &  &  &  &  & 0 \\
\hline
\end{tabular}  
}
}
\end{center}
\end{table}

\noindent Sampling 100000 sets of parameters compatible with the preference information provided by the Dean, the results show that $s_{3}$ and $s_{7}$ are surely the best two students. In fact, looking at the positive and the negative flows separately, that is to how much an alternative outranks all other alternatives and how much an alternative is outranked from all other alternatives, we get that students $s_3$ is preferred to all other students, except to $s_{7}$, with a frequency of at least 68.01\%, while $s_7$ is preferred to all other students, except to $s_3$, with a frequency of at least 71.08\%. Comparing $s_{3}$ and $s_7$, we get instead that $s_{3}$ is preferred to $s_{7}$ with a frequency of 35.7\%, $s_{7}$ is preferred to $s_{3}$ with a frequency of 39.23\% and they result incomparable with a frequency of 25.07\%.

\begin{table}[!h]
\begin{center}
\caption{Results obtained by applying SMAA to the classical PROMETHEE II method}\label{second_table}
\subtable[Preferences in percentage\label{Pref_PROM_II}]{%
\resizebox{0.4\textwidth}{!}{
\begin{tabular}{ccccccccc}
\hline
    &  $\mathbf{s_1}$  &  $\mathbf{s_2}$ &  $\mathbf{s_3}$ &  $\mathbf{s_4}$ &  $\mathbf{s_5}$ &  $\mathbf{s_6}$ &  $\mathbf{s_7}$ &  $\mathbf{s_8}$ \\ 
\hline
 $\mathbf{s_1}$  &        0     & 90,62 & 4,67 & 25,59 & 63,43 & 91,73 & 0     & 47,22 \\
    
 $\mathbf{s_2}$  & 9,38 & 0     & 3,40 & 8,57 & 15,14 & 74,27 & 0     & 20,61 \\
 
 $\mathbf{s_3}$  & 95,33 & 96,6 & 0     & 100   & 88,78 & 93,68 & 53,61 & 87,95 \\
 
 $\mathbf{s_4}$  & 74,41 & 91,43 & 0     & 0     & 70,29 & 89,18 & 28,62 & 59,07 \\
 
 $\mathbf{s_5}$  & 36,57 & 84,86 & 11,22 & 29,71 & 0     & 100   & 1,08 & 37,10 \\
 
 $\mathbf{s_6}$  & 8,27 & 25,73 & 6,32 & 10,82 & 0     & 0     & 0,09 & 8,62 \\
 
 $\mathbf{s_7}$  & 100   & 100   & 46,39 & 71,38 & 98,92 & 99,91 & 0     & 100 \\
 
 $\mathbf{s_8}$  & 52,78 & 79,39 & 12,05 & 40,93 & 62,9 & 91,38 & 0     & 0 \\
\hline
\end{tabular}
}
}
\subtable[First rank acceptability\label{first_rank_acc}]{%
\resizebox{0.45\textwidth}{!}{
\begin{tabular}{ccccccccc}
\hline
    &  $\mathbf{b_{j}^{1}}$  &  $\mathbf{b_{j}^{2}}$ &  $\mathbf{b_{j}^{3}}$ &  $\mathbf{b_{j}^{4}}$ &  $\mathbf{b_{j}^{5}}$ &  $\mathbf{b_{j}^{6}}$ &  $\mathbf{b_{j}^{7}}$ &  $\mathbf{b_{j}^{8}}$ \\ 
\hline
 $\mathbf{s_1}$  &   0     & 0     & 0     & 38,567 & 47,836 & 11,913 & 1,679 & 0,005 \\
    
 $\mathbf{s_2}$  &   0     & 0     & 0     & 5,539 & 13,255 & 5,815 & 57,845 & 17,546 \\
 
 $\mathbf{s_3}$  &   53,612 & 33,045 & 3,211 & 4,033 & 1,164 & 1,786 & 3,149 & 0 \\
 
 $\mathbf{s_4}$  &   0     & 28,625 & 28,896 & 13,881 & 2,778 & 15,222 & 2,468 & 8,13 \\
 
 $\mathbf{s_5}$  &   1,085 & 7,143 & 5,982 & 18,146 & 12,294 & 55,35 & 0     & 0 \\
 
 $\mathbf{s_6}$  &   0     & 0,094 & 2,794 & 3,978 & 1,113 & 3,127 & 19,802 & 69,092 \\
 
 $\mathbf{s_7}$  &   45,303 & 25,978 & 28,719 & 0     & 0     & 0     & 0     & 0 \\
 
 $\mathbf{s_8}$  &   0     & 5,115 & 30,398 & 15,856 & 21,56 & 6,787 & 15,057 & 5,227 \\     
 
\hline
\end{tabular}
}
}

\vspace{0,4cm}
\subtable[Central weight vectors\label{baricentro}]{%
\resizebox{0.3\textwidth}{!}{
\begin{tabular}{cccc}
\hline
\mbox{Students} & $w_1$ & $w_2$ & $w_3$\\
\hline
$s_3$   & 0,61 & 0,22 & 0,17 \\

$s_5$   & 0,02 & 0,42 & 0,56 \\

$s_7$   & 0,12 & 0,62 & 0,26 \\
\hline
\end{tabular}
}
}
\qquad
\subtable[Barycenter\label{bari}]{%
\resizebox{0.2\textwidth}{!}{
\begin{tabular}{ccc}
\hline
   $w_1$ & $w_2$ & $w_3$\\
\hline
   0,25 & 0,42 & 0,33 \\
\hline
\end{tabular}
}
}
\end{center}
\end{table}

\noindent Because the Dean has to give the scholarship to only one of these students, and she wants to be as much cautious as possible in choosing the best one, she would like to have more information about the possible rankings that could be obtained varying the weights assigned to the different subjects. For this reason, we decided to show her also the results obtained by applying SMAA to the classical PROMETHEE II method. We pointed out that, looking at the final ranking of the students obtained with respect to the net flows, student $s_{3}$ is preferred to student $s_{7}$ with a frequency of 53.61\% while $s_7$ is preferred to $s_3$ with a frequency of 43.69\%; looking at Table \ref{first_rank_acc}, $s_7$ will get always one of the first three positions in the final ranking, while $s_3$ will reach the first two positions with a frequency of the 86,65\% but she will reach a position between the fourth and the seventh with a frequency of 10.13\%. \\
At the same time, we observed that only students $s_3$, $s_5$ and $s_7$ can get the first position and this happens if the most important criterion is Mathematics, Literature or Physics respectively (see Table \ref{baricentro}). Looking at Table \ref{bari}, we get also that independently from which alternative reach the first rank, the most important criterion for the DM is Physics followed from Literature and Mathematics.

In the second part of this didactic example, we shall show how to apply the SMAA methodology to the bipolar PROMETHEE methods. \\
As the Dean thinks that scientific subjects are more important than humanistic subjects but, at the same time, she does not want to favour students good in scientific subjects but having lacks in humanistic subjects. Besides, she thinks that Mathematics and Physics are redundant because, generally, a student good in Mathematics is also good in Physics and viceversa. \\
The Dean stated that, locally, student $s_1$ is preferred to student $s_2$ more than student $s_3$ is preferred to student $s_4$, and that student $s_7$ is preferred to student $s_8$ more than student $s_5$ is preferred to student $s_6$. Translating these information with the constraints $Ch^{B}(P^{B}(s_1,s_2),\hat\mu)>Ch^{B}(P^{B}(s_3,s_4),\hat\mu)$ and $Ch^{B}(P^{B}(s_7,s_8),\hat\mu)>Ch^{B}(P^{B}(s_5,s_6),\hat\mu)$ and solving the optimization problem (\ref{first_probl}), we found that $\varepsilon_{1}^{*}<0$, therefore the classical PROMETHEE methods are not able to restore the preference information provided by the Dean. Solving instead the optimization problem (\ref{second_prob}), we found that $\varepsilon_{2}^{*}>0$ and therefore the bipolar PROMETHEE methods are able to restore these preference information.

\begin{table}[!h]
\begin{center}
\caption{Results obtained by the bipolar PROMETHEE I method}
\subtable[Preferences in percentage\label{bip_Pref_PROM_I}]{%
\resizebox{0.4\textwidth}{!}{
\begin{tabular}{ccccccccc}
\hline
    &              $\mathbf{s_1}$  &  $\mathbf{s_2}$ &  $\mathbf{s_3}$ &  $\mathbf{s_4}$ &  $\mathbf{s_5}$ &  $\mathbf{s_6}$ &  $\mathbf{s_7}$ &  $\mathbf{s_8}$ \\ 
\hline
 $\mathbf{s_1}$  &        0 &     15,502 &          0 &          0 &      0,636 &      6,384 &          0 &      0,668 \\

 $\mathbf{s_2}$  &        0 &          0 &          0 &          0 &          0 &      2,844 &          0 &          0 \\

 $\mathbf{s_3}$  &   64,839 &        100 &          0 &     95,552 &     52,853 &     99,897 &      0,242 &     68,591 \\

 $\mathbf{s_4}$  &   13,584 &     86,066 &          0 &          0 &     16,907 &     45,775 &          0 &     30,243 \\

 $\mathbf{s_5}$  &   74,384 &     82,306 &          0 &     33,646 &          0 &     26,957 &          0 &     51,839 \\

 $\mathbf{s_6}$  &        0 &     11,519 &          0 &      0,009 &          0 &          0 &          0 &          0 \\

 $\mathbf{s_7}$  &      100 &        100 &     42,288 &        100 &     99,958 &     99,845 &          0 &        100 \\

 $\mathbf{s_8}$  &   31,673 &     78,226 &          0 &      30,05 &      2,122 &     29,272 &          0 &          0 \\     
 
\hline
\end{tabular}
}
}
\subtable[Incomparability in percentage\label{bip_Incom_PROM_I}]{%
\resizebox{0.4\textwidth}{!}{
\begin{tabular}{ccccccccc}
\hline
    &  $\mathbf{s_1}$  &  $\mathbf{s_2}$ &  $\mathbf{s_3}$ &  $\mathbf{s_4}$ &  $\mathbf{s_5}$ &  $\mathbf{s_6}$ &  $\mathbf{s_7}$ &  $\mathbf{s_8}$ \\
\hline

 $\mathbf{s_1}$  &        0 &     84,498 &     35,161 &     86,416 &      24,98 &     93,616 &          0 &     67,659 \\

 $\mathbf{s_2}$  &     &          0 &          0 &     13,934 &     17,694 &     85,637 &          0 &     21,774 \\

 $\mathbf{s_3}$  &     &           &          0 &      4,448 &     47,147 &      0,103 &      57,47 &     31,409 \\

 $\mathbf{s_4}$  &     &      &      &          0 &     49,447 &     54,216 &          0 &     39,707 \\

 $\mathbf{s_5}$  &      &      &     &      &          0 &     73,043 &      0,042 &     46,039 \\

 $\mathbf{s_6}$  &     &      &     &     &     &          0 &      0,155 &     70,728 \\

 $\mathbf{s_7}$  &        &         &      &          &       &       &          0 &          0 \\

 $\mathbf{s_8}$  &    &     &    &    &     &    &          &          0 \\
\hline
\end{tabular}  
}
}
\end{center}
\end{table}

In Tables \ref{bip_Pref_PROM_I} and \ref{bip_Incom_PROM_I} are shown the results obtained by applying SMAA to the bipolar PROMETHEE I method, while Tables \ref{bip_Pref_PROM_II}, \ref{bip_first_rank_acc} and Table \ref{bip_baricentro} present the results obtained by applying SMAA to the bipolar PROMETHEE II method. By using the bipolar PROMETHEE I method, and therefore considering separately the bipolar positive flow and the bipolar negative flow respectively, we pointed out again, that students $s_3$ and $s_7$ are the best because $s_{3}$ is preferred to all other students, except to $s_7$, with a frequency of at least 52.85\% while $s_7$ is preferred to all other students, except to $s_3$, with a frequency of at least 99.84\%. Comparing $s_3$ and $s_7$, we got that $s_7$ is preferred to $s_3$ with a frequency of 42.28\%, $s_3$ is preferred to $s_7$ with a frequency of 0.242\% and they result incomparable with a frequency of 57.47\%. 

\begin{table}[!h]
\begin{center}
\caption{Results obtained by the bipolar PROMETHEE II method}
\subtable[Preferences in percentage\label{bip_Pref_PROM_II}]{%
\resizebox{0.45\textwidth}{!}{
\begin{tabular}{ccccccccc}
\hline
    &  $\mathbf{s_1}$  &  $\mathbf{s_2}$ &  $\mathbf{s_3}$ &  $\mathbf{s_4}$ &  $\mathbf{s_5}$ &  $\mathbf{s_6}$ &  $\mathbf{s_7}$ &  $\mathbf{s_8}$ \\ 
\hline
 $\mathbf{s_1}$  &         0 &        100 &          0 &     10,193 &      1,998 &     99,267 &          0 &      2,383 \\

 $\mathbf{s_2}$  &         0 &          0 &          0 &          0 &          0 &     32,995 &          0 &          0 \\

 $\mathbf{s_3}$  &       100 &        100 &          0 &        100 &     99,803 &        100 &     11,243 &     99,939 \\

 $\mathbf{s_4}$  &    89,807 &        100 &          0 &          0 &     40,787 &     98,613 &          0 &     50,856 \\

 $\mathbf{s_5}$  &    98,002 &        100 &      0,197 &     59,213 &          0 &        100 &          0 &     86,808 \\

 $\mathbf{s_6}$  &     0,733 &     67,005 &          0 &      1,387 &          0 &          0 &          0 &          0 \\

 $\mathbf{s_7}$  &       100 &        100 &     88,757 &        100 &        100 &        100 &          0 &        100 \\

 $\mathbf{s_8}$  &    97,617 &        100 &      0,061 &     49,144 &     13,192 &        100 &          0 &          0 \\
\hline
\end{tabular}
}
}
\subtable[First rank acceptability\label{bip_first_rank_acc}]{%
\resizebox{0.45\textwidth}{!}{
\begin{tabular}{ccccccccc}
\hline
    &  $\mathbf{b_{j}^{1}}$  &  $\mathbf{b_{j}^{2}}$ &  $\mathbf{b_{j}^{3}}$ &  $\mathbf{b_{j}^{4}}$ &  $\mathbf{b_{j}^{5}}$ &  $\mathbf{b_{j}^{6}}$ &  $\mathbf{b_{j}^{7}}$ &  $\mathbf{b_{j}^{8}}$ \\ 
\hline
 $\mathbf{s_1}$  &   0 &          0 &          0 &      1,569 &     10,938 &     87,258 &      0,235 &          0 \\

 $\mathbf{s_2}$  &        0 &          0 &          0 &          0 &          0 &          0 &     32,995 &     67,005 \\

 $\mathbf{s_3}$  &   11,243 &      88,56 &      0,136 &      0,061 &          0 &          0 &          0 &          0 \\

 $\mathbf{s_4}$  &        0 &          0 &     40,161 &     11,321 &      38,27 &      8,916 &      1,332 &          0 \\

 $\mathbf{s_5}$  &        0 &      0,189 &     54,347 &     36,416 &      7,591 &      1,457 &          0 &          0 \\

 $\mathbf{s_6}$  &        0 &          0 &          0 &          0 &      0,553 &      1,014 &     65,438 &     32,995 \\

 $\mathbf{s_7}$  &   88,757 &     11,243 &          0 &          0 &          0 &          0 &          0 &          0 \\

 $\mathbf{s_8}$  &        0 &      0,008 &      5,356 &     50,633 &     42,648 &      1,355 &          0 &          0 \\     
 
\hline
\end{tabular}
}
}
\end{center}
\end{table}

\begin{table}[!h]
\begin{center}
\subtable[Central weight vectors\label{bip_baricentro}]{%
\resizebox{0.45\textwidth}{!}{
\begin{tabular}{cccccccccc}
\hline
\mbox{Students} & $a_1$ & $a_2$ & $a_3$ & $a_{12}$ & $a_{13}$ & $a_{23}$ & $a^{+}_{1|2}$ & $a^{+}_{1|3}$ & $a^{+}_{2|1}$\\
\hline
$s_3$   & 0,69 &   0,62 &   0,57  &   -0,26 &   -0,32 &   -0,29  & -0,21 &   -0,17 &   -0,08 \\
$s_7$   & 0,56 &   0,69 &   0,64  & -0,29   &   -0,26 &   -0,34  & -0,14 &   -0,16 &   -0,11 \\
\hline
\\
\hline
                & $a^{+}_{2|3}$ & $a^{+}_{3|1}$ & $a^{+}_{3|2}$ & $a^{-}_{1|2}$ & $a^{-}_{1|3}$ & $a^{-}_{2|1}$ & $a^{-}_{2|3}$ & $a^{-}_{3|1}$ & $a^{-}_{3|2}$\\
\hline
$s_3$   & -0,09 &   -0,09 &   -0,10 &   -0,08 &   -0,09&   -0,21 &   -0,10 &   -0,17 &   -0,09 \\
$s_7$   & -0,16 &   -0,05 &   -0,16 &   -0,11 &   -0,05 &   -0,14 &   -0,16 &   -0,16 &   -0,16  \\
\hline
\end{tabular}
}
}
\qquad
\subtable[Barycenter\label{bip_bari}]{%
\resizebox{0.45\textwidth}{!}{
\begin{tabular}{ccccccccc}
\hline
$a_1$ & $a_2$ & $a_3$ & $a_{12}$ & $a_{13}$ & $a_{23}$ & $a^{+}_{1|2}$ & $a^{+}_{1|3}$ & $a^{+}_{2|1}$\\
\hline
 0,625 &   0,655 &   0,605  &   -0,275  &   -0,29 &   -0,315  & -0,175 &   -0,165 &   -0,095 \\
\hline
\\
\hline
                $a^{+}_{2|3}$ & $a^{+}_{3|1}$ & $a^{+}_{3|2}$ & $a^{-}_{1|2}$ & $a^{-}_{1|3}$ & $a^{-}_{2|1}$ & $a^{-}_{2|3}$ & $a^{-}_{3|1}$ & $a^{-}_{3|2}$\\
\hline
-0,125 &   -0,07 &   -0,13 &   -0,095 &   -0,07&    -0,175 &   -0,13 &   -0,165 &   -0,125 \\
\hline
\end{tabular}
}
}
\end{center}
\end{table}

By applying the SMAA methodology to the bipolar PROMETHEE II method, we get that considering the ranking obtained taking into account the bipolar net flow and all sampled sets of parameters, only students $s_3$ and $s_7$ can reach the first position. Besides, they shall fill surely the first two positions but $s_7$ will be the best with a frequency of 88.75\% and, looking at Table \ref{bip_baricentro}, just in case Physics is the most important criterion, followed by Literature and Mathematics.\\
We observe in Table \ref{bip_bari} that, for the DM, the criteria are almost equally important and there is redundancy between all couples of considered criteria.

\section{Conclusions}\label{conclusions}
PROMETHEE methods and Stochastic Multiobjective Acceptability Analysis (SMAA) have been widely applied to deal with several real world problems \cite{behzadian2010promethee,tervonen_figueira}.
In this paper we have proposed to apply SMAA to the classical PROMETHEE methods \cite{Brans_book,Brans84} and to the bipolar PROMETHEE methods \cite{Corrente_annals,Corrente2012a}. The integration of SMAA and PROMETHEE methods permit to get recommendations provided considering several set of parameters compatible with the preference information provided by the Decision Maker (DM), and not considering only one of these sets. Thus the methodology presented in this paper permits to deal effectively with robustness concerns related to the choice of preference parameters in PROMETHEE methods, also in case interaction between criteria and bipolarity of the scales are considered.\\
In the case of PROMETHEE II and of the bipolar PROMETHEE II, one can look at the possible final rankings obtained considering the net flow and the bipolar net flow respectively and at which parameters give to an alternative the best position. \\ 
In the didactic example we have shown that the SMAA methodology can be  applied to PROMETHEE methods, also in cases the criteria are interacting. This brings us to believe that the methodology we are proposing in this paper can be adopted to deal with many complex real-world problems.

\bibliographystyle{plain}
\bibliography{Bib_Integral}

\begin{thebibliography}{10}

\bibitem{angilella_IPMU}
A.~Angilella, S.~Corrente, and S.~Greco.
\newblock {SMAA}-{C}hoquet: {S}tochastic {M}ulticriteria {A}cceptability
  {A}nalysis for the {C}hoquet {I}ntegral.
\newblock In S.~Greco et~al., editor, {\em IPMU 2012, Part IV, CCIS 300}, pages
  248--257. Springer-Verlag Berlin Heidelberg 2012, 2012.

\bibitem{angilella2010non}
S.~Angilella, S.~Greco, and B.~Matarazzo.
\newblock Non-additive robust ordinal regression: A multiple criteria decision
  model based on the {C}hoquet integral.
\newblock {\em European Journal of Operational Research}, 201(1):277--288,
  2010.

\bibitem{behzadian2010promethee}
M.~Behzadian, R.B. Kazemzadeh, A.~Albadvi, and M.~Aghdasi.
\newblock {PROMETHEE}: A comprehensive literature review on methodologies and
  applications.
\newblock {\em European Journal of Operational Research}, 200(1):198--215,
  2010.

\bibitem{Brans_book}
J.P. Brans and B.~Mareschal.
\newblock {PROMETHEE} {M}ethods.
\newblock In J.~Figueira, S.~Greco, and M.~Ehrgott, editors, {\em {Multiple
  Criteria Decision Analysis: State of the Art Surveys}}, pages 163--196.
  Springer, Berlin, 2005.

\bibitem{Brans84}
{J.P.} Brans, B.~Mareschal, and {Ph}. Vincke.
\newblock {PROMETHEE}: a new family of outranking methods in multicriteria
  analysis.
\newblock In {J.P.} Brans, editor, {\em Operational Research, IFORS 84}, pages
  477--490. North Holland, Amsterdam, 1984.

\bibitem{BransVi85}
J.P. Brans and Ph. Vincke.
\newblock A preference ranking organisation method: {T}he {PROMETHEE} method
  for {MCDM}.
\newblock {\em Management Science}, 31(6):647--656, 1985.

\bibitem{choquet1953theory}
G.~Choquet.
\newblock Theory of capacities.
\newblock {\em Ann. Inst. Fourier}, 5(54):131--295, 1953.

\bibitem{Corrente_annals}
S.~Corrente, J.~Figueira, and S.~Greco.
\newblock Dealing with {I}nteraction {B}etween {B}ipolar {M}ultiple {C}riteria
  {P}references in {PROMETHEE M}ethods.
\newblock {\em Annals of Operational Research. Submitted}, 2012.
\newblock \texttt{http://arxiv.org/abs/1211.0507}.

\bibitem{Corrente2012a}
S.~Corrente, J.R. Figueira, and S.~Greco.
\newblock Interaction of criteria and robust ordinal regression in bi-polar
  {PROMETHEE} methods.
\newblock In S.~Greco et~al., editor, {\em IPMU 2012, Part IV, CCIS 300}, pages
  469--479. Springer, Berlin, 2012.

\bibitem{dias2006inferring}
L.C. Dias and V.~Mousseau.
\newblock Inferring electre's veto-related parameters from outranking examples.
\newblock {\em European Journal of Operational Research}, 170(1):172--191,
  2006.

\bibitem{eppe_IPMU}
S.~Eppe and Y.~Se~Smet.
\newblock Studying the {I}mpact of {I}nformation {S}tructure in the {PROMETHEE
  II} {P}reference {E}licitation {P}rocess: {A S}imulation {B}ased {A}pproach.
\newblock In S.~Greco et~al., editor, {\em IPMU 2012, Part IV, CCIS 300}, pages
  383--392. Springer-Verlag Berlin Heidelberg 2012, 2012.

\bibitem{eppe_de_Smet}
S.~Eppe, Y.~Se~Smet, and T.~Stutzle.
\newblock A {B}i-{O}bjective {O}ptimization {M}odel to {E}liciting {D}ecision
  {M}aker's {P}references for the {PROMETHEE II} {M}ethod.
\newblock In R.I. Brafman, F.~Roberts, and A.~Tsoukias, editors, {\em ADT 2011.
  LNCD, vol. 6992}, pages 56--66. Springer Heidelberg 2011, 2011.

\bibitem{FigGreEhr}
J.~Figueira, S.~Greco, and M.~Ehrgott~(eds.).
\newblock {\em Multiple Criteria Decision Analysis: State of the Art Surveys}.
\newblock Springer, Berlin, 2005.

\bibitem{FGR}
J.~Figueira, S.~Greco, and B.~Roy.
\newblock {ELECTRE} methods with interaction between criteria: An extension of
  the concordance index.
\newblock {\em European Journal of Operational Research}, 199(2):478--495,
  2009.

\bibitem{figueira2009building}
J.~Figueira, S.~Greco, and R.~S{\l}owi\'{n}ski.
\newblock Building a set of additive value functions representing a reference
  preorder and intensities of preference: {GRIP} method.
\newblock {\em European Journal of Operational Research}, 195(2):460--486,
  2009.

\bibitem{fujimoto2004new}
K.~Fujimoto.
\newblock New characterizations of k-additivity and k-monotonicity of
  bi-capacities.
\newblock In {\em SCIS-ISIS 2004, 2nd Int. Conf. on Soft Computing and
  Intelligent Systems and 5th Int. Symp. on Advanced Intelligent Systems},
  2004.

\bibitem{Grabisch1996}
M.~Grabisch.
\newblock The application of fuzzy integrals in multicriteria decision making.
\newblock {\em European Journal of Operational Research}, 89:445--456, 1996.

\bibitem{GL1}
M.~Grabisch and C.~Labreuche.
\newblock Bi-capacities-{I}: definition, {M}{\"o}bius transform and
  interaction.
\newblock {\em Fuzzy sets and systems}, 151(2):211--236, 2005.

\bibitem{GL2}
M.~Grabisch and C.~Labreuche.
\newblock Bi-capacities-{II}: the {C}hoquet integral.
\newblock {\em Fuzzy sets and systems}, 151(2):237--259, 2005.

\bibitem{Grabisch_book_greco}
M.~Grabisch and C.~Labreuche.
\newblock Fuzzy measures and integrals in {MCDA}.
\newblock In J.~Figueira, S.~Greco, and M.~Ehrgott, editors, {\em {Multiple
  Criteria Decision Analysis: State of the Art Surveys}}, pages 563--604.
  Springer, Berlin, 2005.

\bibitem{Grabisch2008}
M.~Grabisch and C.~Labreuche.
\newblock A decade of application of the {C}hoquet and {S}ugeno integrals in
  multi-criteria decision aid.
\newblock {\em 4OR}, 6(1):1--44, 2008.

\bibitem{GF2003}
S.~Greco and J.R. Figueira.
\newblock Dealing with interaction between bi-polar multiple criteria
  preferences in outranking methods.
\newblock {\em Research Report 11-2003}, INESC-Coimbra, Portugal, 2003.
\newblock
  \texttt{http://www.uc.pt/en/org/inescc/res\_reports\_docs/research\_reports}.

\bibitem{GKMS2010}
S.~Greco, M.~Kadzi{\'n}ski, V.~Mousseau, and R.~S{\l}owi{\'n}ski.
\newblock {ELECTRE}$^{GKMS}$: Robust ordinal regression for outranking methods.
\newblock {\em European Journal of Operational Research}, 214(1):118--135,
  2011.

\bibitem{GrecoMatarazzoSlowinski02}
S.~Greco, B.~Matarazzo, and R.~Slowinski.
\newblock Bipolar {S}ugeno and {C}hoquet integrals.
\newblock In {\em Proceedings of the EUROFUSE 02 Workshop on Information
  Systems}, pages 191--196, 2002.

\bibitem{greco2008ordinal}
S.~Greco, V.~Mousseau, and R.~S{\l}owi\'{n}ski.
\newblock Ordinal regression revisited: multiple criteria ranking using a set
  of additive value functions.
\newblock {\em European Journal of Operational Research}, 191(2):416--436,
  2008.

\bibitem{GR12}
S.~Greco and F.~Rindone.
\newblock Bipolar fuzzy integrals.
\newblock {\em Fuzzy Sets and Systems}, 2013.
\newblock In press.

\bibitem{greco2010robust}
S.~Greco, R.~S{\l}owi{\'n}ski, J.R. Figueira, and V.~Mousseau.
\newblock Robust ordinal regression.
\newblock In M.~Ehrgott, J.~Figueira, and S.~Greco, editors, {\em Trends in
  Multiple Criteria Decision Analysis}, pages 273--320. Springer, Berlin, 2010.

\bibitem{hokkanen}
J.~Hokkanen, R.~Lahdelma, K.~Miettinen, and P.~Salminen.
\newblock Determining the implementation order of a general plan by using a
  multicriteria method.
\newblock {\em Journal of Multi-Criteria Decision Analysis}, 7(5):273--284,
  1998.

\bibitem{Promethee2010}
M.~Kadzi\'{n}ski, S.~Greco, and R.~S{\l}owi\'{n}ski.
\newblock {Extreme ranking analysis in robust ordinal regression}.
\newblock {\em Omega}, 40(4):488--501, 2012.

\bibitem{Lahdelma}
R.~Lahdelma, J.~Hokkanen, and P.~Salminen.
\newblock {SMAA} - stochastic multiobjective acceptability analysis.
\newblock {\em European Journal of Operational Research}, 106(1):137--143,
  1998.

\bibitem{Lahdelma_S2}
R.~Lahdelma and P.~Salminen.
\newblock {SMAA}-2: Stochastic multicriteria acceptability analysis for group
  decision making.
\newblock {\em Operations Research}, 49(3):444--454, 2001.

\bibitem{mareschal}
B.~Mareschal.
\newblock Weight stability intervals in multicriteria decision aid.
\newblock {\em European Journal of Operational Research}, 33(1):54--64, 1988.

\bibitem{mousseau2003resolving}
V.~Mousseau, J.~Figueira, L.~Dias, C.~Gomes~da Silva, and J.~Climaco.
\newblock Resolving inconsistencies among constraints on the parameters of an
  {MCDA} model.
\newblock {\em European Journal of Operational Research}, 147(1):72--93, 2003.

\bibitem{mousseau1998inferring}
V.~Mousseau and R.~S{\l}owi\'{n}ski.
\newblock Inferring an {ELECTRE TRI} model from assignment examples.
\newblock {\em Journal of Global Optimization}, 12(2):157--174, 1998.

\bibitem{mousseauuser}
V.~Mousseau, R.~S{\l}owi\'{n}ski, and P.~Zielniewicz.
\newblock A user-oriented implementation of the {ELECTRE-TRI} method
  integrating preference elicitation support.
\newblock {\em Computers \& Operations Research}, 27(7-8):757--777, 2000.

\bibitem{smith1984}
R.L Smith.
\newblock Efficient {M}onte {C}arlo procedures for generating points uniformly
  distributed over bounded regions.
\newblock {\em Operations Research}, 32:1296--1308, 1984.

\bibitem{sugeno1974theory}
M.~Sugeno.
\newblock {\em Theory of fuzzy integrals and its applications}.
\newblock {P}h.{D.} {T}hesis, {T}okyo {I}nstitute of {T}echnology, 1974.

\bibitem{sun}
Z.~Sun and M.~Han.
\newblock Multi-criteria decision making based on {PROMETHEE} method.
\newblock In {\em Proceedings of the 2010 International Conference on Computing
  Control and Industrial Engineering, IEEE Computer Society Press}, pages
  416--418, 2010.

\bibitem{tervonen_figueira}
T.~Tervonen and J.~Figueira.
\newblock A survey on stochastic multicriteria acceptability analysis methods.
\newblock {\em Journal of Multi-Criteria Decision Analysis}, 15(1-2):1--14,
  2008.

\bibitem{tervonen_electre_tri}
T.~Tervonen, J.~Figueira, R.~Lahdelma, J.~Almeida~Dias, and P.~Salminen.
\newblock A stochastic method for robustness analysis in sorting problems.
\newblock {\em European Journal of Operational Research}, 192(1):236--242,
  2007.

\bibitem{Tervonen2012}
T.~Tervonen, G.~Van~Valkenhoef, N.~Basturk, and D.~Postmus.
\newblock Hit-and-run enables efficient weight generation for simulation-based
  multiple criteria decision analysis.
\newblock {\em European Journal of Operational Research}, 224:552--559, 2013.

\bibitem{wakker1989additive}
P.P. Wakker.
\newblock {\em Additive representations of preferences: A new foundation of
  decision analysis}, volume~4.
\newblock Springer, 1989.

\bibitem{wolters}
W.T.M. Wolters and B.~Mareschal.
\newblock Novel types of sensitivity analysis for additive {MCDM} methods.
\newblock {\em European Journal of Operational Research}, 81(2):281--290, 1995.

\end{thebibliography}

\end{document}